\documentclass{article}
\usepackage{graphicx}
\usepackage{amsfonts}
\usepackage{amsmath}
\usepackage[dvipsnames]{xcolor}
\usepackage{biblatex}
\usepackage{url}

\title{A New Transformation of the Magic Square of Squares}
\author{Christian Wolird \\ chris.wolird@email.ucr.edu}
\date{July 2023}

\begin{document}

\maketitle
\begin{abstract}
    We show arithmetic triplets of Gaussian squares are in 3-to-1 correspondence with Pythagorean triples thereof. This correspondence would transform a solution to the Magic Square of Squares puzzle into a larger structure of perfect Gaussian squares. In particular, we obtain the backwards result that a puzzle solution would generate non-trivial near-miss solutions in the Gaussian integers. Results are applied to popular near-misses.
\end{abstract}

\section{Triples and Triplets}

Sometimes three perfect squares form an evenly spaced triplet.
\begin{figure}[h]
    \centering
    \includegraphics[scale=0.25]{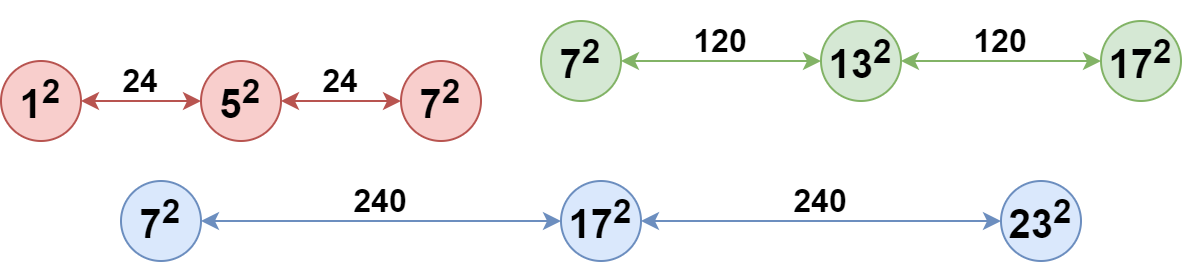}
    \caption{Arithmetic triplets.}
\end{figure}

For brevity, we call such a creature an \textit{arithmetic triplet} (of squares unless specified otherwise). In 1225 Fibonacci published \textit{The Book of Squares} describing, among other things, a connection between these arithmetic triplets and their more popular relative, the Pythagorean triple.

Specifically, there's a 1-to-1 correspondence where the hypotenuse-square of any Pythagorean triple is also the middle-square of an arithmetic triplet (and the other way round). 
\begin{center}
\begin{tabular}{ ccc } 
 $A^2+B^2=C^2$ & $\Rightarrow$ & $(A+B)^2 - C^2 = C^2 - (A-B)^2$ \\
 && \\
 $\big(\frac{L+R}{2}\big)^2+\big(\frac{L-R}{2}\big)^2=C^2$ & $\Leftarrow$ & $L^2-C^2=C^2-R^2$ \\ 
\end{tabular}
\end{center}

\newpage

\begin{figure}[h]
    \centering
    \includegraphics[scale=0.25]{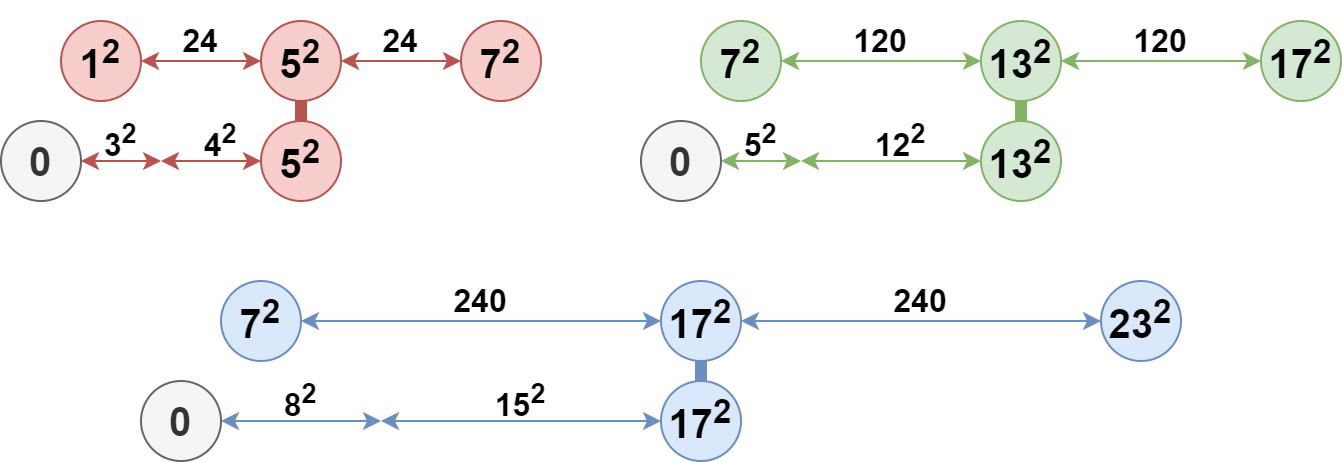}
    \caption{Arithmetic triplet and Pythagorean triple correspondence examples.}
\end{figure}

Because we're concerning ourselves with \textit{integer} solutions 
 here, notice that $\frac{L\pm R}{2}$ must be an integer since $L^2+R^2=2C^2$ means that $L$ and $R$ are both odd or both even.
\begin{figure}[h]
    \centering
    \includegraphics[scale=0.25]{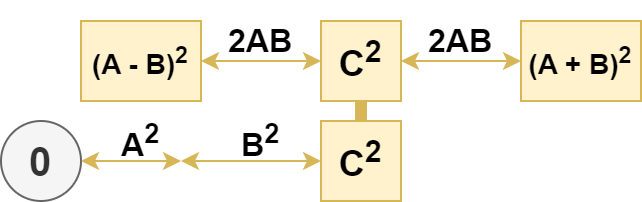}
    \caption{The general arithmetic-Pythagorean correspondence.}
\end{figure}

\section{Unfolding}

Over the integers, this correspondence is like a lawn chair folded flat. To see its real shape, we can unfold in the complex plane.
\begin{figure}[h]
    \centering
    \includegraphics[scale=0.2]{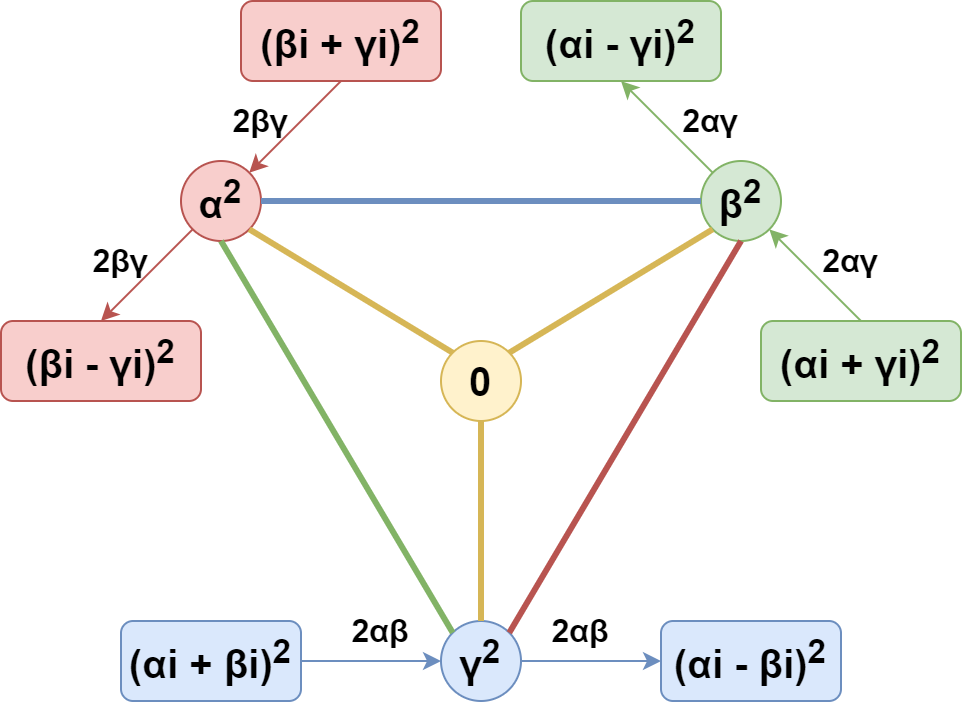}
    \caption{The arithmetic-Pythagorean correspondence ``unfolded" into $\mathbb{Z}[i]$.}
\end{figure}

\newpage

With Pythagorean triples of integers, one square, $C^2$, is doomed to be the hypotenuse and remain forever alone. But with Pythagorean triples of \textit{Gaussian} integers\footnote{Complex numbers with integer parts: $2+i,\ -17+5i,\ 83i,\ $etc.}, the squares can be collected all together.
$${\color{Green}(8-4i)^2}+{\color{blue}(4+7i)^2}={\color{red}(4-i)^2}$$
$$\Downarrow$$
$${\color{red}(4-i)^2}+{\color{Green}(4+8i)^2}+{\color{blue}(7-4i)^2}=0$$

Whereas before, an \textit{integer} Pythagorean triple $({\color{red}A},{\color{Green}B},{\color{blue}C})$ meant a solution to ${\color{red}A^2}+{\color{Green}B^2}={\color{blue}C^2}$, we now think of a \textit{Gaussian} Pythagorean triple $({\color{red}\alpha}, {\color{Green}\beta}, {\color{blue}\gamma})$ as a solution to ${\color{red}\alpha^2}+{\color{Green}\beta^2}+{\color{blue}\gamma^2}=0$, acknowledging no one square doomed to hypotenusity. Thus we create \textit{three} Gaussian arithmetic triplets from any one Gaussian Pythagorean triple by treating each of ${\color{red}\alpha}, {\color{Green}\beta},$ and ${\color{blue}\gamma}$ as the ``hypotenuse" in turn.
 $${\color{red}(4-i)^2}+{\color{Green}(4+8i)^2}+{\color{blue}(7-4i)^2}=0$$
 $$\Downarrow$$
 $${\color{red}(4-i)^2}={\color{Green}(8-4i)^2}+{\color{blue}(4+7i)^2}$$
 $${\color{Green}(4+8i)^2}={\color{red}(1+4i)^2}+{\color{blue}(4+7i)^2}$$
 $${\color{blue}(7+4i)^2}={\color{red}(1+4i)^2}+{\color{Green}(8-4i)^2}$$
 $$\Downarrow$$
 $$(12+3i)^2-{\color{red}(4-i)^2}={\color{red}(4-i)^2}-(4-11i)^2$$
 $$(5+12i)^2-{\color{Green}(4+8i)^2}={\color{Green}(4+8i)^2}-(3+3i)^2$$
 $$(9)^2-{\color{blue}(7+4i)^2}={\color{blue}(7+4i)^2}-(7+8i)^2$$

The arithmetic-Pythagorean correspondence is thus \textit{3}-to-1 in the Gaussians. Any two arithmetic triplets resulting from the same Pythagorean triple, we call \textit{siblings}. The general correspondence can there be visualized (see Figure 4) as a triangle in the complex plane having
\begin{enumerate}
    \item its centroid at zero,
    \item its vertices at perfect integer squares,
    \item and vertices which are themselves the midpoints of three line segments having perfect integer squares for endpoints.
\end{enumerate}

As far as the author could find, the existing treatments of Pythagorean Triples over the Gaussians (such as in [1]) make no mention of these arithmetic triplets.

\newpage

For an instance in the wild, we plot the prior algebraic example:
\begin{figure}[h]
    \centering
    \includegraphics[scale=0.35]{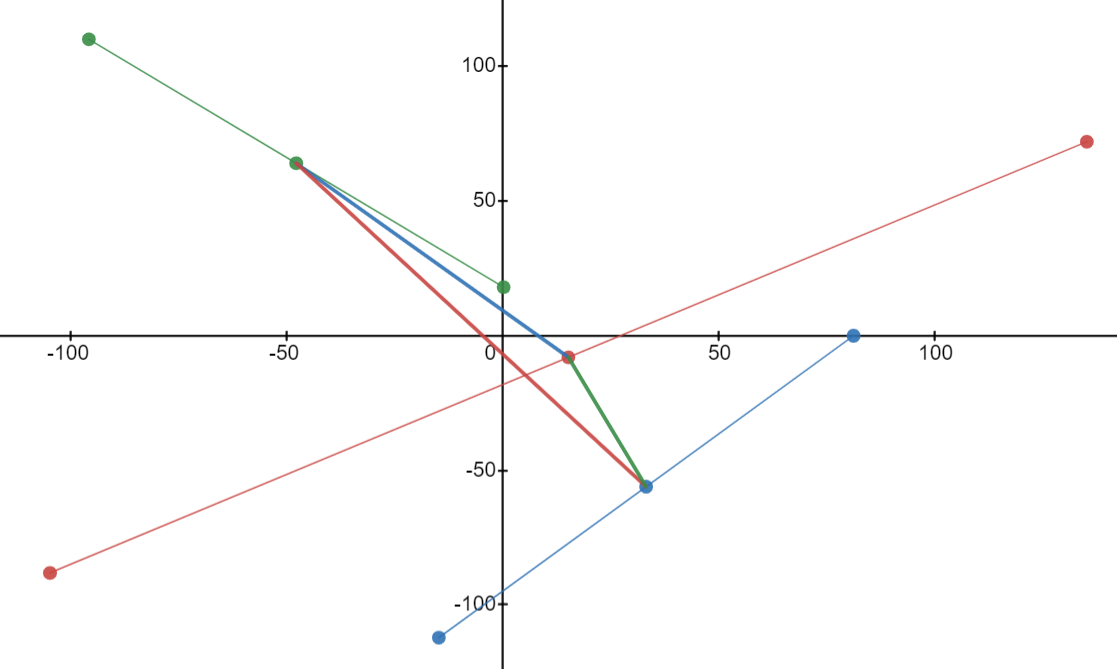}
    \caption{Plot of $(4-i)^2+(4+8i)^2+(7-4i)^2=0$ and the arithmetic triplets it generates.}
\end{figure}

\section{The Magic Square of (Gaussian) Squares}
\begin{figure}[h]
    \centering
    \includegraphics[scale=0.2]{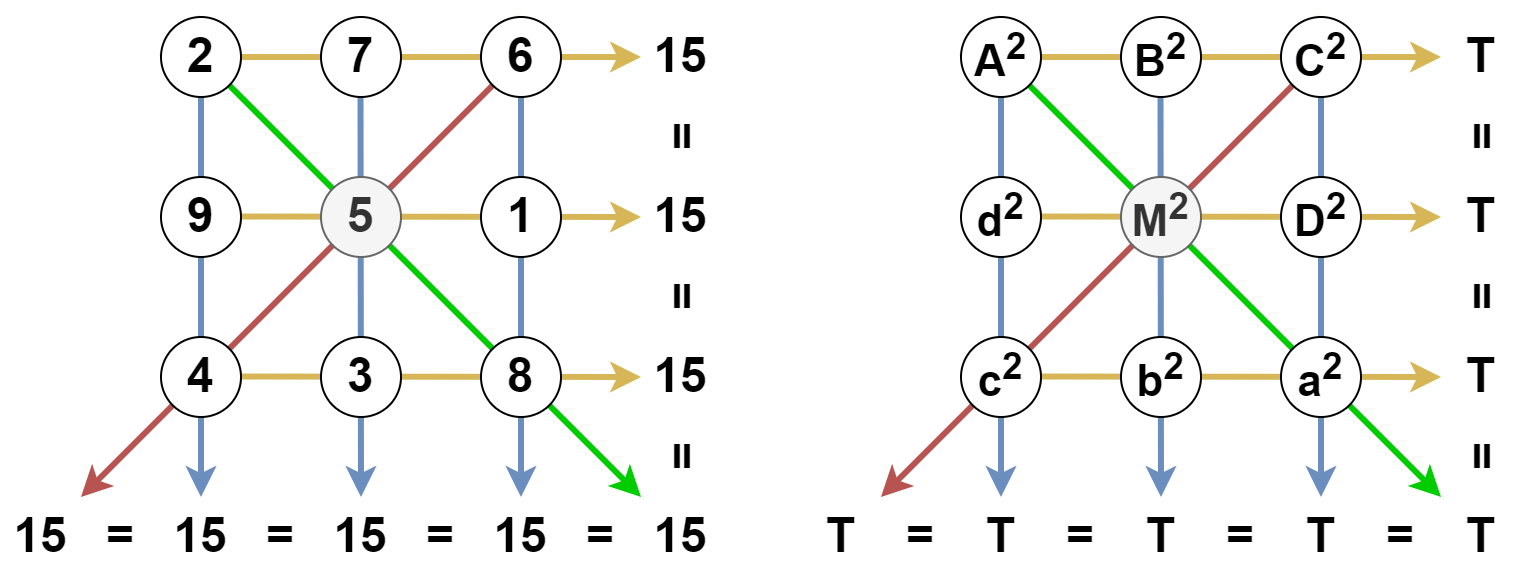}
    \caption{The Lo Shu magic square (left) and the general Magic Square of Squares template (right).}
\end{figure}

The ever-loved and ever-frustrating ``Magic Square of Squares" puzzle is the challenge of arranging nine distinct integer squares in a $3\times3$ grid so that the three numbers in each column, row, and diagonal add to the same total (yet unsolved for both ordinary integers and Gaussians).

\newpage

An important side-effect of these rules is that the magic total must be triple the middle number. By algebraic rearrangement,

$$3T=(A^2+M^2+a^2)+(B^2+M^2+b^2)+(C^2+M^2+c^2)$$
$$=(A^2+B^2+C^2)+(M^2+M^2+M^2)+(a^2+b^2+c^2)=2T+3M^2$$
Or instead, we may offer a more intuitive visual proof:
\begin{figure}[h]
    \centering
    \includegraphics[scale=0.2]{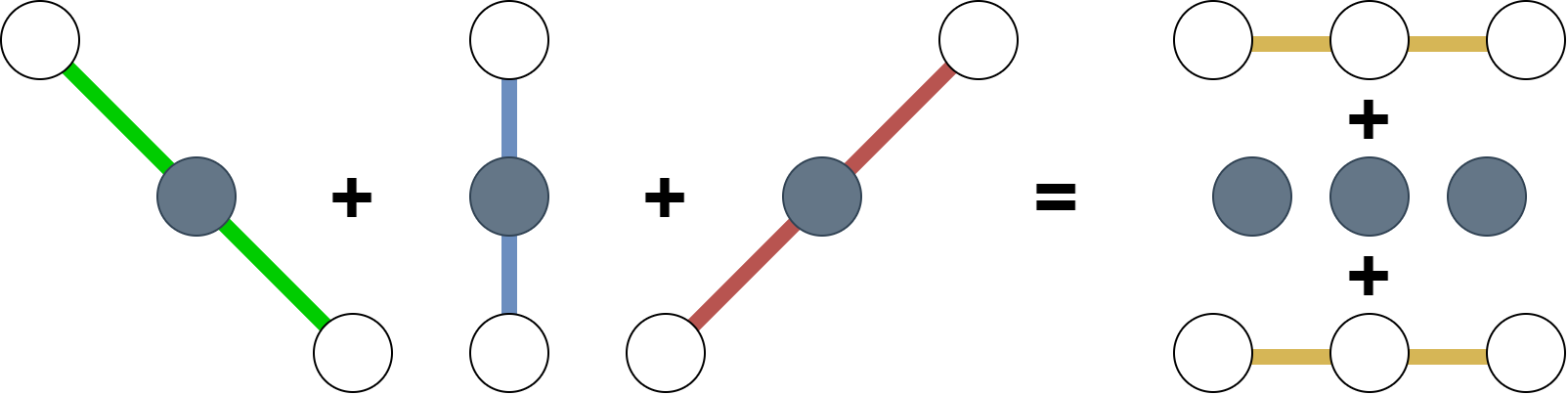}
    \caption{A visual proof that $T=3M^2$.}
\end{figure}

This side-effect in turn means any \textit{central} sum (i.e. the sums containing the center) must form an arithmetic triplet. Taking a diagonal for example,

$$A^2 + M^2 + a^2 = T = 3M^2$$
$$\Downarrow$$
$$A^2 - M^2 = M^2 - a^2$$

All together, this means a magic square of \textit{Gaussian} squares forms a slant $3\times3$ grid\footnote{Called a $3x3$ general arithmetic progression, or GAP, in related articles such as [2].} in the complex plane, containing \textit{8} arithmetic triplets.
\begin{figure}[h]
    \centering
    \includegraphics[scale=0.2]{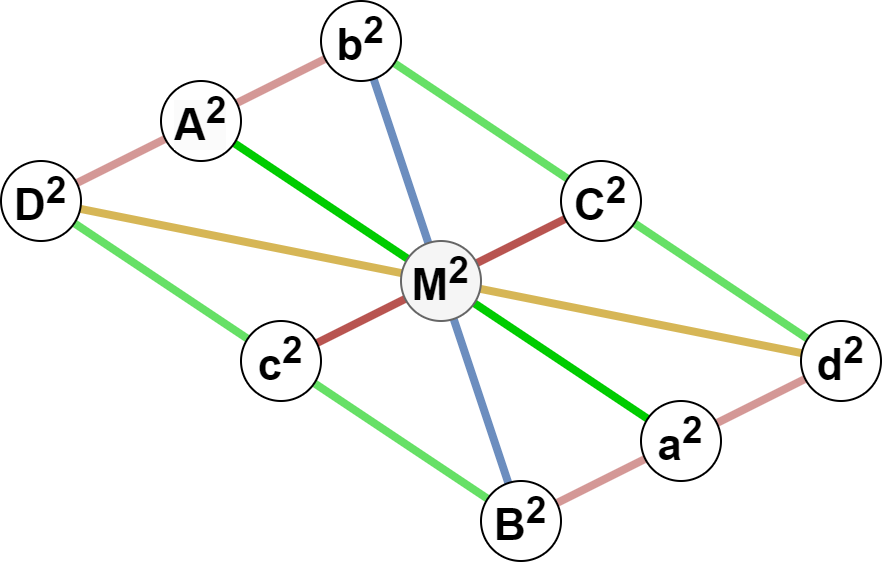}
    \caption{General form of a magic square of Gaussian squares plotted when in the complex plane.}
\end{figure}

\newpage

And from these 8 arithmetic triplets, we may construct another 16 using the 3-to-1 arithmetic-Pythagorean correspondence. Taking the same diagonal for example,
$$A^2-M^2=M^2-a^2$$
$$\Downarrow$$
$$\Big(\frac{A+a}{2}i+M\Big)^2 - \Big(\frac{A-a}{2}\Big)^2 = \Big(\frac{A-a}{2}\Big)^2 - \Big(\frac{A+a}{2}i-M\Big)^2$$
$$\Big(\frac{A-a}{2}i+M\Big)^2 - \Big(\frac{A+a}{2}\Big)^2 = \Big(\frac{A+a}{2}\Big)^2 - \Big(\frac{A-a}{2}i-M\Big)^2$$
Thus any Magic Square of Squares necessarily brings into existence, not 8, but rather 24 Gaussian arithmetic triplets.

To maintain a distinction within these 16 new arithmetic triplets, we call 8 of them the \textit{older} siblings and 8 of them the \textit{younger} siblings. Specifically, the arithmetic triplets having \textit{sums} at their centers, like $\big(\frac{A+a}{2}\big)^2$, are the \textit{older} siblings. And the arithmetic triplets having \textit{differences} at their centers, like $\big(\frac{A-a}{2}\big)^2$, are the \textit{younger} siblings.
\begin{figure}[h]
    \centering
    \includegraphics[scale=0.5]{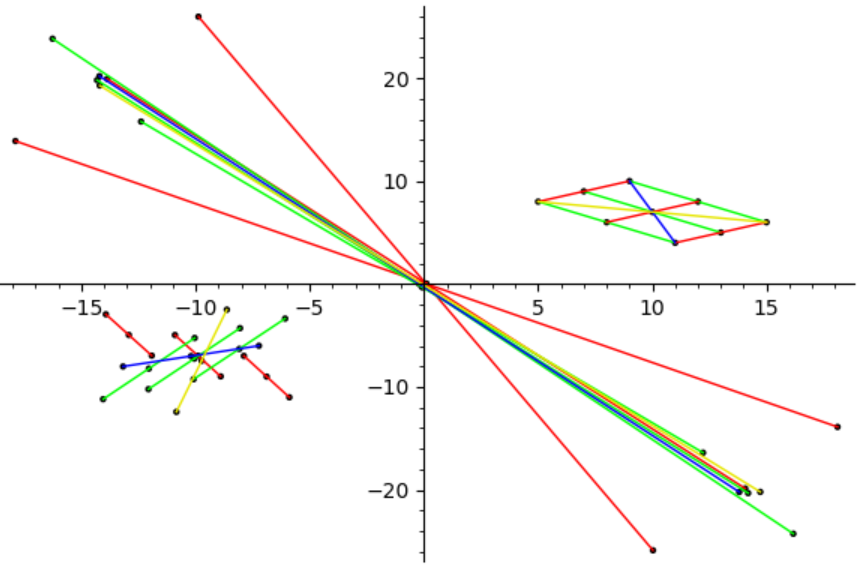}
    \caption{A slant grid (first quadrant) with the older siblings (third quadrant) and younger siblings (second/fourth quadrant) it generates.}
\end{figure}

These siblings can be generated from any $3\times 3$ slant grid in the complex plane. And if the slant grid is composed of perfect squares, then so too will its siblings be. But alas, no $3\times 3$ magic square of Gaussian squares has been found and the visuals presented here are of \textit{non}-perfect squares\footnote{For visual clarity, we are actually displaying a rotation of the older sibling, which would normally also sit in the 1st quadrant. I.e. displaying $-(\frac{A+a}{2})^2=(\frac{A+a}{2i})^2$ instead of $(\frac{A+a}{2})^2$.}.

\newpage

\section{From Solutions to Misses}

There is an abundance of ``near-miss solutions" to the Magic Square of Squares puzzle online[3][4].
\begin{figure}[h]
    \centering
    \includegraphics[scale=0.17]{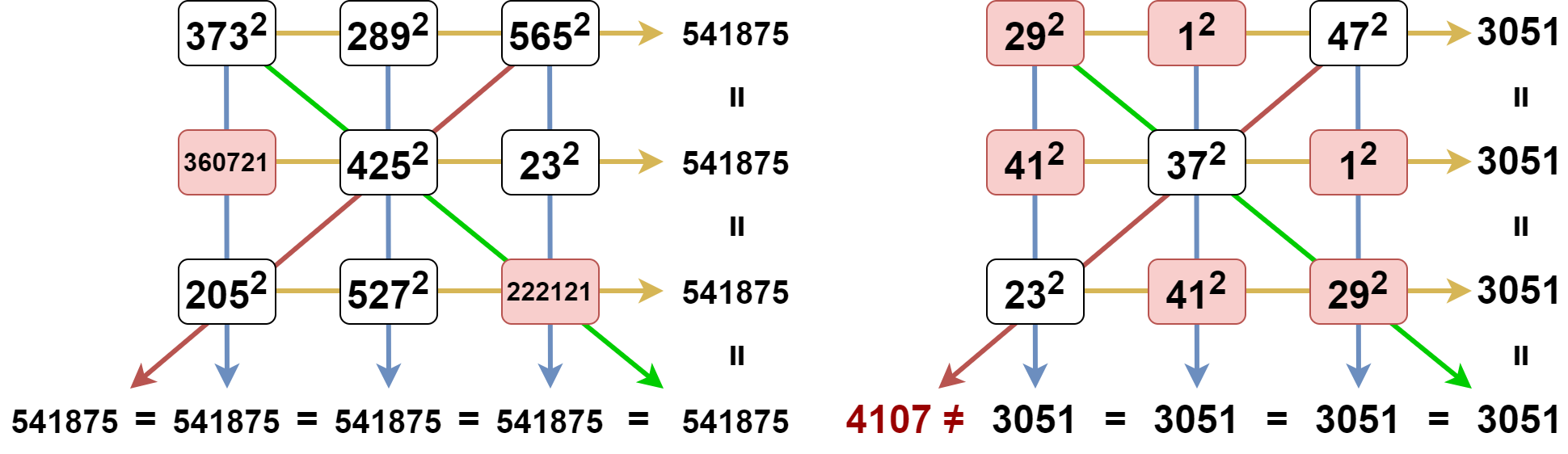}
    \caption{Two well known near-misses: the Bremner square (left) and the Parker square (right).}
\end{figure}

To be clear, by a ``near-miss", we mean a $3\times 3$ grid of numbers that \textit{nearly} meets the constraints of the Magic Square of Squares puzzle, say by having duplicates, a mismatched sum, or a few non-square entries. 

Interestingly, we can find near-misses in the siblings of a true Magic Square of (Gaussian) Squares. Zooming in on the older siblings from the previous section,
\begin{figure}[h]
    \centering
    \includegraphics[scale=0.4]{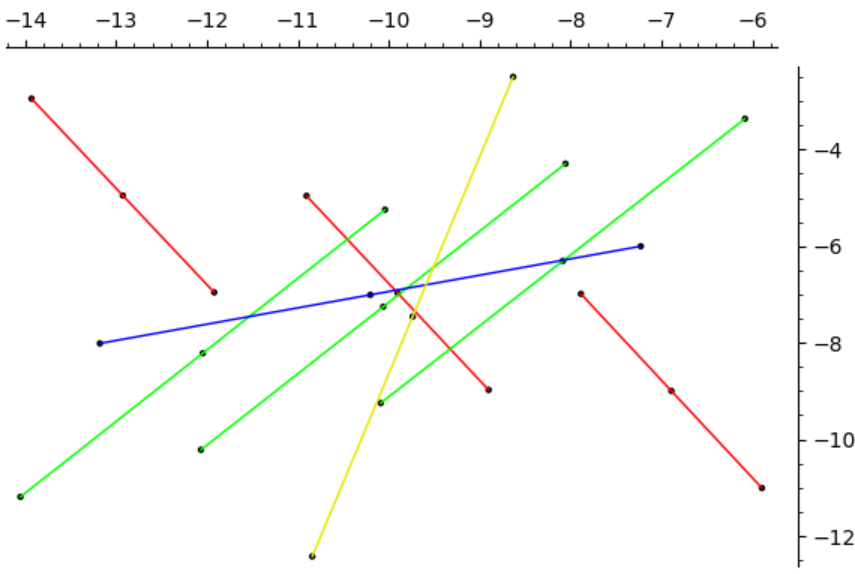}
    \caption{The older siblings of an arbitrary slant grid.}
\end{figure}

\newpage

Six of the arithmetic triplets depicted above seem to form slant grids.
\begin{figure}[h]
    \centering
    \includegraphics[scale=0.4]{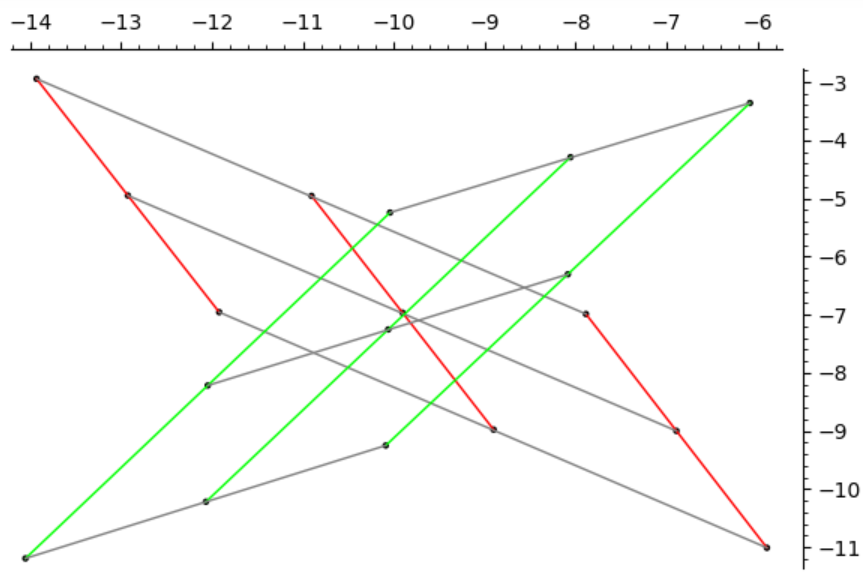}
    \caption{Two near-misses amongst older siblings.}
\end{figure}

So normally, one might solve a puzzle by starting with a near-miss and creating a true solution from it. But our 3-to-1 correspondence indicates that oddly the reverse is true of the Magic Square of Squares. That is, if one finds a solution, there is a good chance it will generate a handful of near-misses. The Magic Square of Squares thus has the mischievous characteristic that it hands you more hints as soon as you've solved it.

But are these really always near-misses? Indeed, from looks alone, it seems we have two new slant grids. Sadly each is slightly kinked. For instance, let's inspect the red pseudo-grid algebraically. The midpoints of the red line segments are
$$\Big(\frac{D+b}{2}\Big)^2,\quad \Big(\frac{c+C}{2}\Big)^2,\quad \text{and}\quad \Big(\frac{B+d}{2}\Big)^2$$
If we had a true slant grid, these three terms ought to form an arithmetic triplet themselves. That is, we ought to have
$$\Big(\frac{D+b}{2}\Big)^2+\Big(\frac{B+d}{2}\Big)^2=2\Big(\frac{C+c}{2}\Big)^2$$
However, the true relationship between these terms is
$$\Big(\frac{D+b}{2}\Big)^2+\Big(\frac{B+d}{2}\Big)^2=\frac{(D^2+B^2)+(d^2+b^2)+2Db+2Bd}{4}$$
$$=\frac{2c^2+2C^2+2Db+2Bd}{4}$$
$$=2\Big(\frac{C+c}{2}\Big)^2+\frac{Db+Bd-2Cc}{2}$$
which differs from a true arithmetic triplet by the error term
$$\frac{Db+Bd-2Cc}{2}$$

To see why $Db+Bd-2Cc$ is nearly zero is algebraically tedious. But the following derivation will do.
$$(Db+Bd-2Cc)(Db+Bd+2Cc)=(Db+Bd)^2-(2C^2)(2c^2)$$
$$=(D^2b^2+2BbDd+B^2d^2)-(B^2+b^2)(D^2+d^2)$$
$$=2BbDd+(D^2b^2+B^2d^2)-(B^2d^2+D^2b^2)-B^2b^2-D^2d^2$$
$$=-(B^2b^2-2BbDd+D^2d^2)=-(Bb-Dd)^2$$
In particular, from this we can rewrite our error term as
$$\frac{Db+Bd-2Cc}{2}=-\frac{(Bb-Dd)^2}{2(Db+Bd+2Cc)}$$

However, this error is only small if our magic square is distanced from the origin. For instance, if we shift the slant grid from the prior figures to overlap with the origin, pseudo-grids are no longer present among the siblings.
\begin{figure}[h]
    \centering
    \includegraphics[scale=0.5]{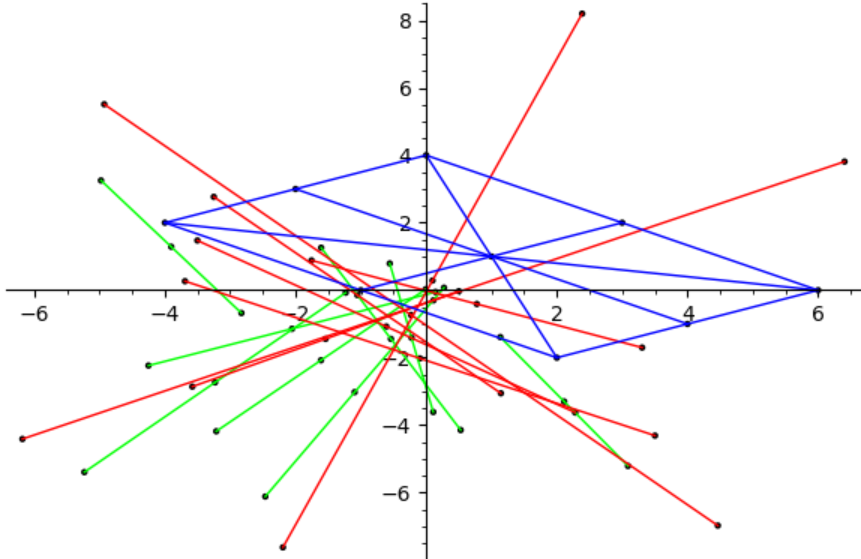}
    \caption{A slant grid covering the origin (blue) along with its older/younger siblings (green/red) amongst which, there are no near-misses.}
\end{figure}

\newpage

\section{Back to Reality}
Any normal Magic Square of Squares sits on the real number line. And the real number line sits in the complex plane. So we can even generate the \textit{complex} siblings of a \textit{real} Magic Square of Squares. As no true solution has yet been found, we show here the siblings of two near-misses, the Bremner square and the Parker square (see Figure 10).
\begin{figure}[h]
    \centering
    \includegraphics[scale=0.42]{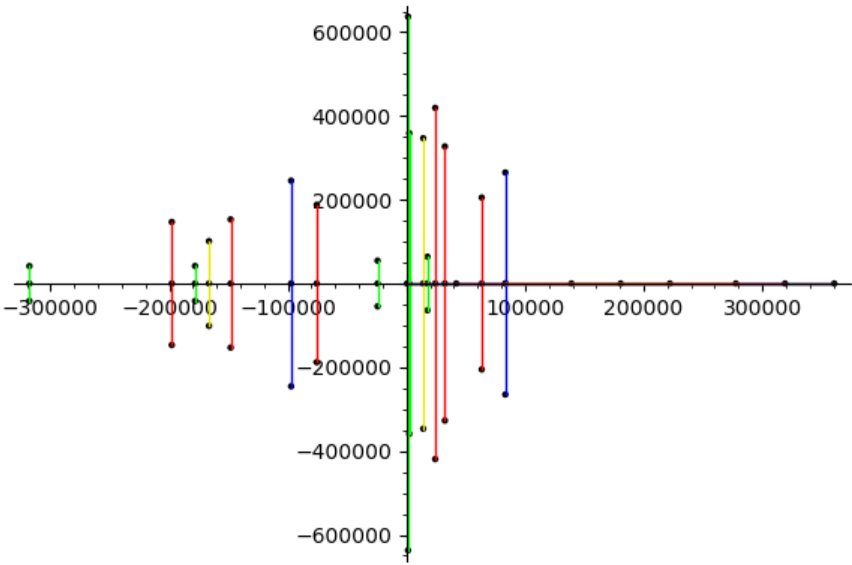}
    \caption{Siblings of the Bremner square.}
\end{figure}

\begin{figure}[h]
    \centering
    \includegraphics[scale=0.42]{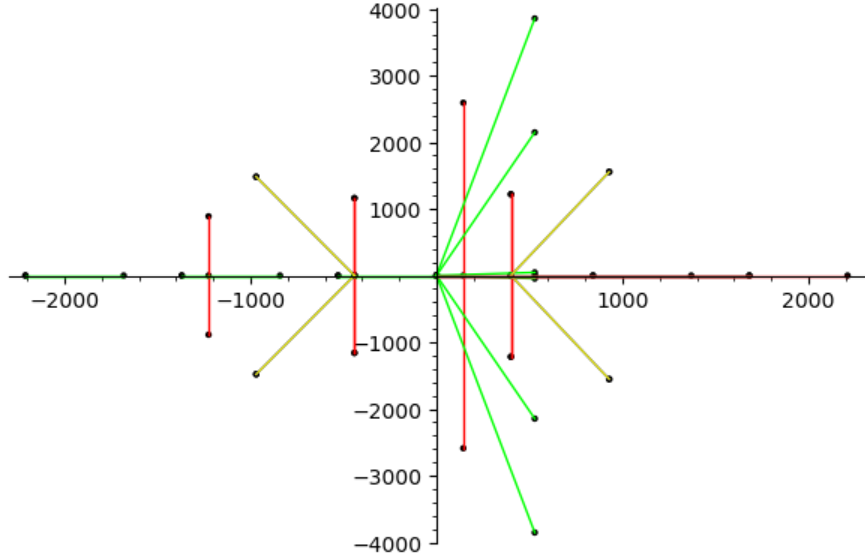}
    \caption{Siblings of the Parker square.}
\end{figure}

\newpage

The interesting bit here is that the personality of each near-miss carries over to its siblings. The Bremner Square has perfect sums and thus its siblings are made of perfect arithmetic triplets. But the Bremner Square has two non-square entries and so its siblings consist only partially of perfect squares.

The Parker Square, on the other hand, has perfectly square entries and thus its siblings are made entirely of perfect Gaussian squares. However, the Parker Square fails to have perfect sums and thus its siblings have a few kinks. And the Parker square fails to have totally distinct entries and thus so do its siblings.

So do these siblings tell us anything about the existence of the Magic Square of Squares? Not that the author sees directly. The author simply wished to make known here the rather spirited personality which the Magic Square of Squares seems to take on in the Gaussians. And the author would be glad to see more serious analyses of the puzzle over the Gaussians in the future.

\end{document}